\newcommand{\CR}{\hbox{{$\cal R$}}}
\newcommand{\CF}{\hbox{{$\cal F$}}}
\newcommand{\CC}{\hbox{{$\cal C$}}}
\newcommand{\CQ}{\hbox{{$\cal Q$}}}

\newcommand{\cg}{\mathfrak{g}}

\newcommand{\R}{\mathbb{R}}
\newcommand{\C}{\mathbb{C}}

\newcommand{\h}{{\scriptstyle\frac{1}{2}}}
\newcommand{\extd}{{\rm d}}
\newcommand{\del}{\partial}
\newcommand{\isom}{{\cong}}
\newcommand{\eps}{{\epsilon}}
\newcommand{\tens}{\mathop{\otimes}}
\newcommand{\la}{{\triangleright}}
\newcommand{\End}{{\rm End}}
\newcommand{\id}{{\rm id}}
\newcommand{\<}{\langle}

\newcommand{\eqn}[2]{\begin{equation}#2\label{#1}\end{equation}}

\newcommand{\rcross}{{\triangleright\!\!\!<}}
\newcommand{\codcross}{{\blacktriangleright\!\!\blacktriangleleft}}
\newcommand{\bicross}{{\blacktriangleright\!\!\!\triangleleft}}

\renewcommand{\o}{{{}_{\scriptscriptstyle(1)}}}
\renewcommand{\t}{{{}_{\scriptscriptstyle(2)}}}
\renewcommand{\th}{{{}_{\scriptscriptstyle(3)}}}
\newcommand{\bo}{{{}^{\scriptscriptstyle(\bar 1)}}}
\newcommand{\bt}{{{}^{\scriptscriptstyle(\bar 2)}}}

\documentclass[10pt]{article}
\usepackage{amssymb,amsmath,epsfig}
\newtheorem{lemma}{Lemma}
\newtheorem{propos}[lemma]{Proposition}

\newtheorem{theorem}[lemma]{Theorem}

\newtheorem{corol}[lemma]{Corollary}

\begin{document}%\baselineskip 22pt

{\ }\qquad \hskip 4.3in
\vspace{.2in}

\begin{center} {\LARGE Noncommutative Differential Geometry and
Twisting of Quantum Groups}
\\ \baselineskip 13pt{\ }\\
{\ }\\ Shahn Majid \\ {\ }\\ School of Mathematical Sciences,
Queen Mary and Westfield College\\ University of London, Mile End
Rd, London E1 4NS\footnote{Reader and Royal Society University
Research Fellow}

\end{center}
%\vspace{10pt}

\begin{quote}%\baselineskip 13pt
\noindent{\bf Abstract}
We outline the recent classification of differential structures for
all main classes of quantum groups. We also outline the algebraic
notion of `quantum manifold' and `quantum Riemannian manifold' based on
quantum group principal bundles, a formulation that works over general
unital algebras.
 \end{quote}

\section{Introduction}

There have been many attempts in the last decades to arrive at a
theory of noncommutative geometry applicable to `coordinate'
algebras that are not necessarily commutative, notably that of A.
Connes coming out of abstract $C^*$-algebra theory in the light of
the Gelfand-Naimark and Serre-Swan theorems. One has tools such as
cyclic cohomology and examples such as the noncommutative torus
and other foliation $C^*$-algebras. Another `bottom up' approach,
which we outline, is based on the idea that the theory should be
guided by the inclusion of the large vein of `naturally occuring'
examples, the coordinate algebras of the quantum groups $U_q(g)$
in particular, and Hopf algebras in general, whose validity for
several branches of mathematics has already been established. This
is similar to the key role that Lie groups played in the
development of modern differential geometry. Much progress has
been made in recent years and there is by now (at least at the
algebraic level) a more or less clear formulation of `quantum
manifold' suggested by this approach. After being validated on the
q-deformation examples such as quantum groups, quantum homogeneous
spaces etc, one can eventually apply the theory quite broadly to a
wide range of unital algebras. The approach will be algebraic,
although not incompatible with $C^*$ completions at a later stage.
In particular, as a bonus, one can apply the theory to
finite-dimensional algebra, i.e. to discrete classical and quantum
systems.

An outline of the paper is the following. We start with the lowest
level structure which (in our approach) is the choice of
differential structure. This is the topic of Section~2 where we
outline the recently achieved more or less complete classification
results. In Section~3 we develop the notion of `quantum
manifold'\cite{Ma:rie} based on noncommutative frame bundles with
quantum group fibre. Usual ideas of `sheaf theory' and `local
trivialisations' do not work in this setting, but from
\cite{BrzMa:gau} one has global algebraic replacements. There is also
an appropriate notion of automorphism or diffeomorphism quantum
groups\cite{Ma:dif}.

\section{Quantum differential forms}

Let $M$ be a unital algebra, which we consider as playing the role of
`co-ordinates' in algebraic geometry, except that we do not require
the algebra to be commutative. The appropriate notion of cotangent
space or differential 1-forms in this case is

1. $\Omega^1$ an $M$-bimodule

2. $\extd:M\to \Omega^1$ a linear map obeying the Leibniz rule
$\extd(ab)=a\extd b+(\extd a)b$ for all $a,b\in M$.

3. The map $M\tens M\to \Omega^1$, $a\tens b\mapsto a\extd b$ is
surjective.

Differential structures are not
unique even classically, and even more non-unique in the quantum
case. There is, however, one universal example of which others are
quotients. This is
\eqn{univ}{ \Omega^1_{\rm univ}=\ker\cdot\subset M\tens M,\quad \extd
a=a\tens 1-1\tens a.}
This is common to more or less all approaches to noncommutative geometry.

The main feature here is that, in usual algebraic geometry, the
multiplication of forms $\Omega^1$ by `functions' $M$ is the same
from the left or from the right. However, if $a\extd b=(\extd b)a$
then by axiom 2. we have $\extd(ab-ba)=0$, i.e. we cannot naturally
suppose this when $M$ is noncommutative.  We say that a
differential calculus is noncommutative or `quantum'  if the left
and right multiplication of forms by functions do not coincide.

When $M$ has a Hopf algebra structure with coproduct $\Delta:M\to
M\tens M$ and counit $\eps:M\to k$ ($k$ the ground field), we say
that $\Omega^1$ is {\em bicovariant} if\cite{Wor:dif}

4. $\Omega^1$ is a bicomodule with coactions $\Delta_L:\Omega^1\to
M\tens\Omega^1,\Delta_R:\Omega^1\to \Omega^1\tens M$ bimodule maps (with
the tensor product bimodule structure on the target spaces, where
$M$ is a bimodule by left and right multiplication).

5. $\extd$ is a bicomodule map with the left and right regular
coactions on $M$ provided by $\Delta$.

A morphism of calculi means a bimodule and bicomodule map forming a
commuting triangle with the respective $\extd$ maps. One
says\cite{Ma:cla} that a calculus is {\em coirreducible} if it has no
proper quotients. Whereas the translation-invariant calculus is unique
classically, in the quantum group case we have at least complete
classification results in terms of representation theory\cite{Ma:cla}.
The dimension of a calculus is that of its space of (say)
left-invariant 1-forms, which can be viewed as generating the rest of
the calculus as a right $M$-module. Similarly with left and right
interchanged.

We note also that in the bicovariant case there is a natural
extension\cite{Wor:dif} from $\Omega^1$ to $\Omega^n$ with $\extd^2=0$.
This is defined as the tensor algebra over $M$ generated by $\Omega^1$
modulo relations defined by a braiding which acts by a simple transposition
on left-invariant and right-invariant forms. Other extensions are also
possible and in general the differential structure
can be specified order by order. Given the extension, one has a quantum DeRahm
cohomology defined in the usual way as closed forms modulo exact ones. Apart
from cohomology one can also start to do `$U(1)$' gauge theory with trivial
bundles, where a gauge field is just a differential form $\alpha\in \Omega^1$
and its curvature is $F=\extd\alpha+\alpha\wedge\alpha$, etc. A gauge transform
is \eqn{trivialgauge}{ \alpha^\gamma=\gamma^{-1}\alpha\gamma
+\gamma^{-1}\extd \gamma,\quad F^\gamma=\gamma^{-1} F\gamma} for any
invertible `function' $\gamma\in M$, and so on. One can define then the
space of flat connections as those with $F=0$ modulo gauge
transformation.  This gives two examples of `geometric' invariants
which work therefore for general algebras equipped with differential
structure.

\subsection{$M=k[x]$}

For polynomials in one variable the
coirreducible calculi have the
form\cite{Ma:fie}
\[ \Omega^1=k_\lambda[x],\quad  \extd
f(x)=\frac{f(x+\lambda)-f(x)}{\lambda},\quad f(x)\cdot
g(\lambda,x)=f(x+\lambda)g(\lambda,x),\quad g(\lambda,x)\cdot
f(x)=g(\lambda,x)f(x)\]
for functions $f$ and one-forms $g$.  Here $k_\lambda$ is a field
extension of the form $k[\lambda]$ modulo $m(\lambda)=0$ and $m$ is an
irreducible monic polynomial.  The dimension of the calculus is the
order of the field extension or the degree of $m$.

For example, the calculi on $\C[x]$ are classified by $\lambda_0\in
\C$ (here $m(\lambda)=\lambda-\lambda_0$) and one has
\[ \Omega^1=\extd x\C[x],\quad \extd f=\extd
x\frac{f(x+\lambda_0)-f(x)}{\lambda_0},\quad x\extd x=(\extd
x)x+\lambda_0.\] We see that the Newtonian case $\lambda_0=0$ is
only one special point in the moduli space of quantum differential
calculi. But if Newton had not supposed that differentials and
forms commute he would have had no need to take this limit. What
one finds with noncommutative geometry is that there is no need to
take this limit at all.  It is also
interesting that the most important field extension in physics,
$\R\subset\C$, can be viewed noncommutative-geometrically with
complex functions $\C[x]$ the quantum 1-forms on the algebra of
real functions $\R[x]$. There is nontrivial quantum DeRahm cohomology
in this case.

\subsection{$M=\C[G]$}

For the coordinate algebra of a finite group $G$ (for convenience
we work over $\C$)  the coirreducible calculi correspond to nontrivial
conjugacy classes $\CC\subset G$ and have the form
\[ \Omega^1=\CC\cdot\C[G],\quad \extd f
=\sum_{g\in\CC}g\cdot (L_g(f)-f),\quad f\cdot g=g\cdot L_g(f)\]
where $L_g(f)=f(g\cdot)$ is the translate of $f\in \C[G]$. The dimension
of the calculus is the order of the conjugacy class.

For the coordinate algebra $\C[G]$ of a Lie group with Lie algebra
$\cg$ the coirreducible calculi  correspond to maximal ideals in
$\ker\eps$ stable under the adjoint coaction. Or in a natural
reformulation\cite{Ma:cla} in terms of quantum tangent spaces the
correspondence is with irreducible $Ad$-invariant subspaces of the
enveloping algebra $\ker\eps\subset U(\cg)$ which are stable under the
coaction $\Delta_L=\Delta-\id\tens 1$ of $U(\cg)$. For example $\cg$
itself defines the standard translation-invariant calculus and this is
coirreducible when $\cg$ is semisimple.

\subsection{$M=\C G$}

For the group algebra of a nonAbelian finite group $G$, we definitely
need the machinery of noncommutative geometry since $M$ itself is
noncommutative. We regard these group algebras `up side down' as if
coordinates, i.e. we describe the geometry of the noncommutative space
$\hat G$ in some sense. The above definitions make sense and
differential structures abound.  The coirreducible calculi correspond
to pairs $(V,\rho,\lambda)$ where $(V,\rho)$ is a nontrivial
irreducible representation and $\lambda\in V/\C$\cite{Ma:cla}. They
have the form \[ \Omega^1=V\cdot \C G,\quad \extd  g
=((\rho(g)-1)\lambda)\cdot g,\quad g\cdot v=(\rho(g)v)\cdot g\] where
$g\in G$ is regarded as a `function'. The dimension of the calculus is
that of $V$. The minimum assumption for merely a differential calculus
is that $\lambda$ should be cyclic.

For $M=U(\cg)$ (the Kirillov-Kostant quantisation of $\cg^*$) one has a
similar construction for any irreducible representation $V$ of the Lie
algebra $\cg$ and choice of ray $\lambda$ in it. Then \[\Omega^1=V\cdot
U(\cg),\quad \extd \xi=\rho(\xi)\lambda,\quad \xi\cdot
v=\rho(\xi)v+v\cdot\xi\] where $\xi\in\cg$ is regarded as a `function'.
The dimension is again that of $V$.

For example, let $\cg=b_+$ be the 2-dimensional Lie algebra with
$[x,t]=x$. Let $V$ be the 2-dimensional representation with matrix and
ray vector
\[ \rho(t)=\begin{pmatrix}0&0\\ 0&1\end{pmatrix},\quad
\rho(x)=\begin{pmatrix}0&1\\ 0&0\end{pmatrix},\quad
\lambda=\begin{pmatrix}0\\ 1\end{pmatrix}.\]
Then $\extd t=\lambda$ and $\extd x$ are the usual basis of $V$ and obey
\[ [t,\extd x]=[x,\extd x]=0,\quad [t,\extd t]=\extd t,\quad
[x,\extd t]=\extd x.\] Replacing $x$ trivially by a vector $x_i$,
$i=1,2,3$ one obtains similarly a natural candidate for noncommutative
Minkowski space along with its differential structure. It has measurable
astronomical predictions\cite{AmeMa:wav}.

This covers the classical objects or their duals viewed `up side down'
as noncommutative spaces. For a finite-group bicrossproduct
$\C[M]\bicross \C G$ the classification is a mixture of the two cases
above and is given in \cite{BegMa:dif}. The Lie version remains to be
worked out in detail. The important example of the Planck scale quantum
group $\C[x]\bicross\C[p]$, however, is a twisting by a cocycle of its
classical limit $\C[\R\rcross\R]$ and is therefore covered by a later
subsection.

\subsection{Proofs} The above cases are all sufficiently elementary
that they can be easily worked out using the following simple
observations known essentially (in some form or other) since
\cite{Wor:dif}. We suppose for convenience that $H$ has invertible
antipode.

1. $\ker\eps\subset M$ is an object in the braided category of left
crossed $M$-modules (i.e. modules over the quantum double $D(M)$ in the
finite-dimensional case)  by multiplication and the left adjoint
coaction.

2. The isomorphism $M\tens M\isom M\tens M$ given by $a\tens b\mapsto
(\Delta a)b$ restricts to an isomorphism $\Omega^1_{\rm Univ}\isom
\ker\eps \tens M$ of bimodules and of bicomodules, where the right hand
side is a right (co)module by the (co)product of $M$ and a left
(co)module by the tensor product of the (co)action on $\ker\eps$ and
the (co)product of $M$.

This implies that every other bicovariant $\Omega^1$ is of the form
$\Omega^1\isom \Omega_0\tens M$ where $\Omega_0$ is a quotient object
of $\ker\eps$ in the category of crossed $M$-modules. I.e. the calculi
correspond to ideals in $\ker\eps$ stable under the adjoint coaction.
Given $\Omega^1$ the space $\Omega_0$ is given by the right-invariant
differentials. In categorical terms the braided category of bicovariant
$M$-bimodules as featuring above (i.e. bimodules which are also
bicomodules with structure maps being bimodule maps) can be identified
with that of crossed $M$-modules, under which the Hopf module for the
universal calculus corresponds to $\ker\eps$.

When this is combined with the notion of coirreducibility and with the
Peter-Weyl decomposition of an appropriate type for $\ker\eps$, one
obtains the classification results above. These latter steps have been
introduced by the author\cite{Ma:cla} (before that one found only
sporadic examples of calculi on particular quantum groups, usually
close to the unique classical calculus.)

We also note that for any finite-dimensional bicovariant calculus the
map $\extd:M\to \Omega_0\tens M$ can be viewed as a `partial
derivative' $\del_x:M\to M$ for each $x\in \Omega_0^*$. The space
$\Omega_0^*$ is called the invariant `quantum tangent space' and is
often more important than the 1-forms in applications. These $\del_x$
are not derivations but together form a braided derivation in the
braided category of $M$-crossed modules (there is a braiding as $x\in
\Omega_0^*$ passes $a\in M$) \cite{Ma:cla}.

\subsection{Cotriangular quantum groups and twisting of calculi}

We recall\cite{Ma:book} that if $M$ is a quantum group and $\chi:M\tens
M\to k$ a cocycle in the sense
\[ \chi(b\o\tens c\o)\chi(a\tens b\t c\t)=\chi(a\o\tens b\o)
\chi(a\t b\t\tens c),\quad \chi(1\tens a)=\eps(a),\quad\forall a,b,c\in M\]
then there is a `twisted' quantum
group $M^\chi$ with product \[ a\cdot_\chi b=\chi(a\o\tens b\o)a\t b\t
\chi^{-1}(a\th\tens b\th)\] and unchanged unit, counit and coproduct.
Here $\chi^{-1}$ is the inverse in $(M\tens M)^*$, which we assume, and
$\Delta a=a\o\tens a\t$, etc., is a notation.

\begin{theorem}\cite{MaOec:twi} The bicovariant differentials
$\Omega^1(M^\chi)$ are in 1-1 correspondence with those of $M$.
\end{theorem}
In fact the entire exterior algbera in the bicovariant
case is known to be a super-Hopf algebra (Brzezinski's theorem) and
that of $M^\chi$ is the twist of that of $M$ when $\chi$ is trivially
extended to a cocycle on the latter. The more direct proof involves the
following:

\begin{theorem}\cite{MaOec:twi} There is an equivalence $\CF$ of braided
monoidal categories from left $M$-crossed modules to left
$M^\chi$-crossed modules given by the functor
\[ \CF(V,\la,\Delta_L)=(V,\la^\chi,\Delta_L),\quad a\la^\chi v
=\chi(a\o\tens v\bo)(a\t\la v\bt)\bt\chi^{-1}((a\t\la v\bt)
\bo\tens a\th),\quad \forall a\in M,\ v\in V,\]
where $\la$ denotes the action and $\Delta_Lv=v\bo\tens v\bt$ is a
notation. There is an associated natural transformation
\[ c_{V,W}: \CF(V)\tens \CF(W)\isom \CF(V\tens W),\quad c_{V,W}(v\tens w)
=\chi(v\bo\tens w\bo)v\bt\tens w\bt.\]
\end{theorem}

As a corollary we deduce by Tannaka-Krein reconstruction arguments:

\begin{corol}\cite{MaOec:twi} When $M$ is finite-dimensional the dual
of the Drinfeld double $D(M^\chi)$ is isomorphic to a twist  of the
dual of the Drinfeld double $D(M)$ by $\chi^{-1}$ viewed on $D(M)\tens
D(M)$. \end{corol}

Here we use the theorem\cite{Ma:tan} that an equivalence of comodule
categories respecting the forgetful functor corresponds to a twist of
the underlying quantum groups. The corollary itself can then be
verified directly at an algebraic level once the required (nontrivial)
isomorphism has been found in this way. Of course one can state it also
in terms of Drinfeld's coproduct twists.

Starting with a classical (commutative) Hopf algebra such a twist
yields a cotriangular one and (from recent work of Etingof and
Gelaki\cite{EtiGel}) every finite-dimensional cotriangular Hopf algebra
in the (co)semisimple case over $k$ algebraically closed is of this
form. Hence the differential calculus in this case reduces by the above
theorem to the classification in the classical cases considered in
previous sections. There are many other instances where an important
quantum group is a twisting of another -- the theorem provides its
differential calculus from that of the other.

\subsection{Factorisable quantum groups $M=\C_q[G]$}

Finally, we come to the standard quantum groups $\C_q[G]$ dual to the
Drinfeld-Jimbo $U_q(\cg)$.  Here\cite{Ma:cla} the coirreducible calculi
are essentially provided by nontrivial finite-dimensional irreducible
right comodules $V$ of the quantum group (i.e. essentially by the
irreducible representations of the Lie algebra) and have the form
\eqn{diffqG}{ \Omega^1=\End(V)\cdot\C_q[G],\quad \extd a
=\rho_+(a\o)\circ \rho_-(Sa\t)\cdot a\th-\id\cdot a,\quad
a\cdot\phi=\rho_+(a\o)\circ\phi\circ \rho_-(S a\t)\cdot a\th }
for all $\phi\in\End(V)$, where
\[ \rho_+(a)v=v\bo\CR(a\tens v\bt),\quad \rho_-(a)v
=v\bo\CR^{-1}(v\bt\tens a)\]
and $\CR:\C_q[G]\tens \C_q[G]\to \C$ is the dual-quasitriangular
structure\cite{Ma:book}. $S$ denotes the antipode. The construction
works for any dual-quasitriangular Hopf algebra with factorisable $\CR$
(the minimum one needs for a differential calculus is that
$\CQ_\rho(a)=\rho_+(a\o)\rho_-(Sa\t)$ is surjective) and gives a
classification of calculi if $M$ has in addition the Peter-Weyl
property that $M=\oplus_V V\tens V^*$ as left and right $M$-comodules.
This is the result in \cite[Thm. 4.3]{Ma:cla} cast in a purely comodule
form. Or in the original formulation in terms of quantum tangent spaces
the correspondence is given in more familiar terms with representations
$\rho$ of the quantum enveloping algebra and
\[ \rho_+(a)=(a\tens\rho)(\CR),\quad \rho_-(a)=(\rho\tens a)(\CR^{-1}),\]
where $\CR$ is the quasitriangular structure or universal R-matrix of
$U_q(\cg)$. For finite-dimensional representations only a finite number
of terms in its powerseries contribute here, i.e. there are no
convergence problems.

The factorisability etc. hold formally for $\C_q[G]$ so that although
one has one natural calculus for each irreducible representation there
are also some `shadows' or technical variants allowed according to the
precise formulation of the relevant quantum groups and their duality
(this is more a deficit in the technical definitions than anything
else). The latter aspect has been subsequently clarified in
\cite{BauSch}\cite{HecSch} following our basic result in
\cite{Ma:cla}.

For the sake of a sketch of the proof of the basic result assume that
$M$ is strictly factorisable dual-quasitriangular and has Peter-Weyl
decomposition in terms of irreducible representations $(V,\rho)$ of a
suitable dual Hopf algebra $U$.  Classifying the quotient $M$-crossed
modules of $\ker\eps$ is equivalent essentially to classifying the
subobjects of $\ker\eps\subset U$ as $U$-crossed modules (the quantum
tangent spaces). When $U$ is strictly factorisable its quantum double
$D(U)$ is isomorphic to $U\codcross U$ which, as an algebra, is a
tensor product (the coproduct is twisted). Hence $U$-crossed modules
are equivalent to two $U$-modules. Next, under the isomorphism of
linear spaces $U\isom M$ assumed under strict factorisability, this is
the same as classifying subobjects of $\ker\eps\subset M$. The
$U$-crossed module structure on this final $\ker\eps\subset M$ under
this chain of reasoning is simply evaluation against $M$ coacting
independently from the left and the right (viewed from the left via the
antipode). This is just the action with respect to which the assumed
Peter-Weyl decomposition $M=\oplus_V \End(V)$ is a decomposition into
irreducibles as $V$ runs over the irreducible representations of $U$.
One may make a similar proof working only with $M$-crossed modules
and $M$-comodules throughout and the corresponding comodule Peter-Weyl
decomposition.

\begin{propos}
\cite{Ma:cla} The quantum tangent spaces $\Omega_0^*=V^*\tens V$ for
the above differential calculi on $\C_q[G]$ are braided-Lie algebras
in the sense of \cite{Ma:lie}. The action of basis element $f^i\tens
e_i$ is
\[ \del_{x^i{}_j} (a)=\< x^i{}_j,a\o\>a\t,\quad x^i{}_j
=\CQ(\rho^i{}_j)-\delta^i{}_j\]
where $\CQ:\C_q[G]\to U_q(\cg)$ is defined by $\CR_{21}\CR$ and
$\rho^i{}_j\in \C_q[G]$  are the matrix elements in the representation
$V$ with basis $\{e_i\}$ and dual basis $\{f^i\}$.
\end{propos}

Recall that the dual of any quantum group acts on the quantum group by
the `coregular representation' in the manner shown, in our case by the
$x^i{}_j$. These and $1$ together form a braided-Lie algebra. This is
described by a system of axioms in any braided category including a
pentagonal `braided-Jacobi' identity. Moreover, such objects have
braided enveloping algebras which, for usual Lie algebras $\cg$,
returns a homogenized version of $U(\cg)$. In our case it returns a
quadratic and braided version of $U_q(\cg)$, i.e. this solved (some
years ago\cite{Ma:lie}) the Lie problem for such quantum groups. The
above gives its geometric interpretation.

For the example of $\C_q[SU_2]$ there is basically one bicovariant
calculus for each spin $j$ with dimension $(2j+1)^2$. The lowest
corresponds to the 4-dimensional braided-Lie algebra $gl_{q,2}$ spanned
by
\[ x=\begin{pmatrix}q^H-1& q^{-\h}(q-q^{-1})q^{H\over 2}X_-\\
q^{-\h}(q-q^{-1})X_+q^{H\over 2}&
q^{-H}-1+q^{-1}(q-q^{-1})^2X_+X_-\end{pmatrix} \]
in the usual notations for $U_q(su_2)$. This braided-Lie algebra is
irreducible for generic $q$ but as $q\to 1$ it degenerates into
$su_2\oplus u(1)$. The partial derivatives degenerate into the usual
invariant vector fields on $SU_2$ and an additional 2nd order operator
which turns out to be the Casimir or Laplacian.

\subsection{Discrete manifolds}

To close with one non-quantum group example, consider
any actual manifold with a finite good cover $\{U_i\}_{i\in I}$.
Instead of building geometric invariants on a manifold and
studying them modulo diffeomorphisms we can use the methods above
to first pass to the skeleton of the manifold defined by its open
set structure and do differential geometry directly on this
indexing set $I$. Thus we take $M=\C[I]$ which just means
collections $\{f_i\in \C\}$. The universal $\Omega^1$ is just
matrices $\{f_{ij}\}$ vanishing on the diagonal. We use the
intersection data for the open sets to set some of these to zero.
Similarly for higher forms. Thus\cite{BrzMa:dif}
\[\Omega^1=\{f_{ij}|\ U_i\cap U_j\ne \emptyset\},\quad
\Omega^2=\{f_{ijk}|\ U_i\cap U_j\cap U_k\ne\emptyset\}\]
\[(\extd f)_{ij}=f_i-f_j,\quad (\extd f)_{ijk}
=f_{ij}-f_{ik}+f_{jk}\]
and so on. Then one has that the quantum cohomology is just the
additive Cech cohomology of the original manifold. Similarly,
 one has that the zero curvature gauge
fields modulo gauge transformations recovers again the first Cech
cohomology, but now in a multiplicative form.

\section{Bundles and connections}

The next layer of differential geometry is bundles, connections, etc.
Usually in physics one needs only the local picture with trivial
bundles in each open set -- but for a general noncommutative algebra
$M$ there may be no reasonable `open sets' and one has therefore to
develop the global picture from the start. This is needed for example
to describe the frame bundle of a topologically nontrivial `manifold'.
It also turns out to be rather easier to do the gauge theory beyond the
`$U(1)$' case (i.e. with a nontrivial quantum structure group and
nonuniversal calculus on it) if one takes the global point of view,
even if the bundle itself is trivial. We take a Hopf algebra $H$ in the
role of `functions' on the structure group of the bundle. To keep
things simple we concentrate on the universal differential calculus but
it is important that the general case is also covered by making
suitable quotients.  Recall that a classical bundle has a free action
of a group on the total space $P$ and a local triviality property. In
our algebraic terms we need \cite{BrzMa:gau}:

1. An algebra $P$ and a coaction $\Delta_R:P\to P\tens H$ of the
quantum group $H$ such that the fixed subalgebra is $M$,
\eqn{fixed}{ M=P^H=\{p\in P|\ \Delta_R p=p\tens 1\}.}
We assume that $P$ is flat as an $M$-module.

2.  The sequence
\eqn{exactness}{ 0\to P(\Omega^1M)P\to\Omega^1P{\buildrel{\rm ver}
\over\longrightarrow} P\tens \ker\eps\to 0}
is exact, where ${\rm ver}=(\cdot\tens\id)\Delta_R$.

The map ${\rm ver}$ plays the role of generator of the vertical vector
fields corresponding classically to the action of the group (for each
element of $H^*$ it maps $\Omega^1P\to P$ like a vector field).
Exactness on the left says that the one-forms $P(\Omega^1M)P$ lifted
from the base are exactly the ones annihilated by the vertical vector
fields.  In the universal calculus case this can be formulated as a
Hopf-Galois extension, a condition arising in other contexts in Hopf
algebra theory also. The differential geometric picture is more
powerful and includes general calculi when we use the right-handed
version of $\Omega_0$ in place of $\ker\eps$.

One can then define a connection as an equivariant splitting
\eqn{connection}{ \Omega^1P=P(\Omega^1 M)P\oplus {\rm complement}}
i.e. an equivariant projection $\Pi$ on $\Omega^1P$. One can
show\cite{BrzMa:gau} the required analogue of the usual theory, i.e.
that such a projection corresponds to a connection form such that
\eqn{conform}{\omega:\ker\eps\to\Omega^1P,\quad {\rm ver}\circ
\omega=1\tens\id}
where $\omega$ intertwines with the adjoint coaction of $H$ on
itself. Finally,  if $V$ is a
vector space on which $H$ coacts then we define the associated
`bundles' $E^*=(P\tens V)^H$ and $E=\hom^H(V,P)$, the space of
intertwiners. The two bundles should be viewed geometrically as
`sections' in classical geometry of bundles associated to $V$ and
$V^*$. Given a suitable (so-called strong) connection one has a covariant
derivative\cite{BrzMa:gau}
\eqn{covderiv}{ D_\omega:E\to \Omega^1 M\tens_M E,\quad D_\omega
=(\id-\Pi)\circ \extd.}

All of this can be checked out for the example of
the $q$-monopole bundle over the
$q$-sphere\cite{BrzMa:gau}. Recall that classically the inclusion
$U(1)\subset SU_2$ in the diagonal has coset space $S^2$ and
defines the $U(1)$ bundle over the sphere on which the monopole
lives. In our case the coordinate algebra of $U(1)$
is the polynomials $\C[g,g^{-1}]$ and the classical inclusion
becomes the projection
\[\pi:\C_q[SU_2]\to \C[g,g^{-1}],\quad
\pi\left(\begin{matrix}
a&b\\ c&d
\end{matrix}\right)=\left(\begin{matrix}
g&0\\ 0&g^{-1}
\end{matrix}\right)\]
Its induced coaction $\Delta_R=(\id\tens\pi)\Delta$ is by the
degree defined as the number of $a,c$ minus the number of $b,d$ in
an expression. The quantum sphere $\C_q[S^2]$ is the fixed
subalgebra i.e. the degree zero part. Explicitly, it is generated
by $b_3=ad$, $b_+=cd$, $b_-=ab$ with $q$-commutativity relations
\[ b_\pm b_3=q^{\pm 2}b_3b_\pm+(1-q^{\pm 2})b_\pm,
\quad q^{2}b_-b_+=q^{-2}b_+b_-+(q-q^{-1})(b_3-1)\]
and the sphere equation $b_3^2=b_3+qb_-b_+$. When $q\to 1$ we can
write $b_\pm=\pm(x\pm\imath y)$, $b_3=z+\h$ and the sphere
equation becomes $x^2+y^2+z^2=\frac{1}{4}$ while the others become
that $x,y,z$ commute. One may verify that we have a quantum bundle
in the sense above and that there is a connection
$\omega(g-1)=d\extd a-qb\extd c$ which, as $q\to 1$, becomes the
usual Dirac monopole constructed algebraically. If we take $V=k$
with coaction $1\mapsto 1\tens g^n$, the sections of the
associated vector bundles $E_n$ for each charge $n$ are just the
degree $n$ parts of $\C_q[SU_2]$. The associated covariant
derivative acts on these.

This example demonstrates compatibility with the more traditional
$C^*$-algebra approach of A. Connes\cite{Con:geo} and others.
Traditionally a vector bundle over any algebra is defined as a finitely
generated projective module. However, there was no notion of quantum
principal bundle before quantum groups.

\begin{propos}\cite{HajMa:pro} The
associated bundles $E_n$ for the $q$-monopole bundle
are finitely generated projective modules, i.e. there exist
\[ e_n\in M_{|n|+1}(\C_q[S^2]),\quad e_n^2=e_n,\quad E_n
=e_n\C_q[S^2]^{|n|+1}.\]
The covariant derivative for the monopole has the form $e_n\extd e_n$.
The classes $[e_n]$ are elements of the noncommutative $K$-theory
$K_0(\C_q[S^2])$ and have nontrivial duality pairing with cyclic
cohomology, hence the $q$-monopole bundle is nontrivial.
\end{propos}

The potential applications of quantum group gauge theory hardly
need to be elaborated.  For example, for a classical manifold
\eqn{pi1}{ \left\{ {{{\rm Flat\ connections\ on}\ G-{\rm bundle}}
\atop{\rm modulo\ gauge}}\right\}\isom \hom(\pi_1,G)/G}
using the holonomy. One can view this as a functor from groups to sets
and the homotopy group $\pi_1$ as more or less the representing object
in the category of groups. The same idea with quantum group gauge
theory essentially defines a homotopy quantum group $\pi_1(M)$ for any
algebra $M$ as more or less the representing object of the functor that
assigns to a quantum group $H$ the set of zero-curvature gauge fields
with this quantum structure group.  This goes somewhat beyond vector
bundles and $K$-theory alone.  Although in principle defined, this idea
has yet to be developed in a computable form.

Finally we mention that one needs to make a slight generalisation of
the above to include other noncommutative examples of interest.  In
fact (and a little unexpectedly) the general theory above can be
developed with only a coalgebra rather than a Hopf algebra $H$.  Or
dually it means only an algebra $A$ in place of the enveloping algebra
of a Lie algebra. This was achieved more recently, in
\cite{BrzMa:coa}\cite{BrzMa:geo}, and allows us to include the full
2-parameter quantum spheres as well as (in principle)  all known
$q$-deformed symmetric spaces.  This setting of gauge theory based on
inclusions of algebras could perhaps be viewed as an algebraic analogue of the
notion of `paragroup' in the theory of operator algebras. Also, in a
different direction, one may do the quantum group gauge theory in any
braided monoidal category at the level of braids and
tangles\cite{Ma:diag} so that one has braided group gauge theory and in
principle gauge theory for quasiassociative algebras such
as\cite{AlbMa:qua} the octonions.

\section{Quantum soldering forms and metrics}

We are finally ready to take the plunge and offer at least a first
definition of a `quantum manifold'.   The approach we take is basically
that of the existence of a bein or, in global terms, a soldering form.
The first step is to define a generalised frame bundle or {\em frame
resolution} of our algebra $M$ as\cite{Ma:rie}

1. A quantum principal bundle
$(P,H,\Delta_R)$ over $M$.

2. A comodule $V$ and an equivariant
`soldering form' $\theta:V\to P\Omega^1M\subset\Omega^1P$ such
that the induced map
\eqn{frame}{ E^*\to \Omega^1M,\quad p\tens v\mapsto p\theta(v)}
is an isomorphism.

What this does is to express the cotangent
bundle as associated to a principal one. Other tensors are then
similarly associated, for example vector fields are
$E\isom\Omega^{-1}M$. Of course, all of this has to be done with
suitable choices of differential calculi on $M,P,H$ whereas we
have been focusing for simplicity on the universal calculi. There
are some technicalities here but more or less the same definitions
work in general. The working definition\cite{Ma:rie} of a {\em
quantum manifold} is simply this data
$(M,\Omega^1,P,H,\Delta_R,V,\theta)$. The definition works in that
one has analogues of many usual results. For example, a connection
$\omega$ on the frame bundle induces a covariant derivative
$D_\omega$ on the associated bundle $E^*$ which maps over under
the soldering isomorphism to a covariant derivative
\eqn{nabla}{ \nabla:\Omega^1M\to \Omega^1M\tens_M\Omega^1M.}
Its torsion is defined as corresponding similarly to
$\bar D_\omega\theta$, where we use a suitable (right-handed)
version of the covariant derivative.

Defining a Riemannian structure can be done in a `self-dual' manner as
follows. Given a framing, a `generalised metric' isomorphism
$\Omega^{-1}M\isom\Omega^1M$ between vector fields is equivalent
to\cite{Ma:rie}

3. Another framing $\theta^*:V^*\to(\Omega^1M)P$, which we
call the {\em coframing}, this time with $V^*$.

The associated quantum metric is
\eqn{quametric}{ g= \theta^*(f^a)\theta(e_a)\in\Omega^1M
\tens_M\Omega^1M}
where $\{e_a\}$ is a basis of $V$ and $\{f^a\}$ is a dual basis.

Now, this self-dual formulation of `metric' as framing and coframing is
symmetric between the two. One could regard the coframing as the
framing and vice versa. From our original point of view its torsion
tensor corresponding to $D_\omega\theta^*$ is some other tensor, which
we call the {\em cotorsion tensor}. A natural proposal for a
generalised Levi-Civita connection on a quantum Riemannian manifold is
therefore\cite{Ma:rie}

4. A connection $\omega$ such that the torsion and cotorsion
tensors both vanish.

There is a corresponding covariant derivative $\nabla$. The Riemannian
curvature of course corresponds to the curvature of $\omega$, which is
$\extd\omega+\omega\wedge\omega$, via the soldering form. I would not
say that the Ricci tensor and Einstein tensor are understood abstractly
enough in this formalism but of course one can just write down the
relevant contractions and proceed blindly.

\begin{theorem}\cite{Ma:rie}
Every quantum group $M$ has a framing by $H=M$, $P=M\tens M$,
$V=\ker\eps$ and $\theta$ induced from the quantum group Maurer-Cartan
form $e(v)=Sv\o\tens v\t$. Likewise for all $M$ equipped with a
bicovariant differential calculus, with $V=\Omega_0$.
\end{theorem}

In this construction one builds the framing from a $V$-bein $e$
inducing the isomorphism $\Omega^1(M)=M\tens\Omega_0$ as in Section~2
(in a right-handed setting). Moreover, for quantum groups such as
$\C_q[SU_2]$ there is an Ad-invariant non-degenerate braided Killing
form\cite{Ma:lie} on the underlying braided-Lie algebra,
 which provides a coframing from a framing -- so that quantum groups
such as $\C_q[SU_2]$ with the corresponding bicovariant differential
calculi are quantum Riemannian manifolds in the required sense. The
existence of a generalised Levi-Civita connection in such cases remains
open and may require one to go beyond strong connections.

At least with the universal calculus every quantum homogeneous space is
a quantum manifold too. That includes quantum spheres, quantum planes
etc. In fact, there is a notion of comeasuring or quantum automorphism
bialgebra\cite{Ma:dif} for practically any algebra $M$ and when this
has an antipode (which typically requires some form of completion) one
can write $M$ as a quantum homogeneous space. So almost any algebra $M$
is more or less a quantum manifold for some principal bundle (at least
rather formally).  This is analogous to the idea that any classical
manifold is, rather formally, a homogeneous space of diffeomorphisms
modulo diffeomorphisms fixing a base point.

Finally, to get the physical meaning of the cotorsion tensor and
other novel ideas coming out of this noncommutative Riemannian geometry,
let us consider the semiclassical limit. What we find is that
noncommutative geometry forces us to slightly generalise
conventional Riemannian geometry itself \cite{Ma:rie}:

1. We should allow any group $G$ in the `frame bundle', hence
the more general concept of a `frame resolution'
$(P,G,V,\theta_\mu^a)$ or {generalised manifold}.

2. The {generalised metric} $g_{\mu\nu}
=\theta^*_\mu{}^a\theta_{\nu a}$
corresponding to a coframing $\theta^*_{\mu}{}^a$ is nondegenerate
but need not be symmetric.

3. The {generalised Levi-Civita} connection defined as having
vanishing torsion and vanishing cotorsion respects the metric only
in a skew sense
\eqn{genlev}{ \nabla_\mu g_{\nu\rho}-\nabla_\nu g_{\mu\rho}=0}
and need not be uniquely determined.

This generalisation of Riemannian geometry includes special cases
of symplectic geometry, where the generalised metric is totally
antisymmetric.  It is also remarkable that metrics with antisymmetric
part are exactly what are needed in string theory to establish
T-duality. In summary, one has on the table a general
noncommutative Riemannian geometry to play with. It can be applied to
a variety of algebras far removed from conventional geometry. Some finite
dimensional examples will be presented elsewhere.

\subsection*{Acknowledgements} It is a pleasure to thank the
organisers of the LMS symposium in Durham for a thoroughly
enjoyable conference.

%\bibliographystyle{unsrt}
%\bibliography{biblio}

\end{document}